%% file: bounded.tex
\documentstyle[amsfonts,12pt,thmsb]{article}
\input psfig
\input imsmark.tex

\def\e2piipq{e^{2\pi i p/q}}
\input{tcilatex}

\begin{document}

\title{Bounded Hyperbolic Components of Quadratic Rational Maps}
\author{Adam Lawrence Epstein}
\date{September 1, 1997}

\maketitle

\begin{abstract}
Let ${\cal H}$ be a hyperbolic component of quadratic rational maps
possessing two distinct attracting cycles. We show that ${\cal H}$ has
compact closure in moduli space if and only if neither attractor is a fixed
point.
\end{abstract}

\SBIMSMark{1997/9}{September, 1997}{}
\tableofcontents

\section{Introduction}

From the perspective of dynamics, the simplest rational maps are {\em %
hyperbolic}: every critical point tends under iteration to some attracting
periodic cycle. Such maps constitute an open and conjecturally dense set in
parameter space \cite{mss}, whose components are referred to as {\em %
hyperbolic components}. Maps in the same component are quasiconformally
conjugate near the Julia set, and thus have essentially identical dynamics
if critical orbit relations are ignored.

The family $P_{c}(z)=z^{2}+c$ of quadratic polynomials contains one
unbounded component, namely ${\Bbb C}-M$ where 
\[
M=\{c:\;J(P_{c})\mbox{ is connected}\,\} 
\]
is the much-studied Mandelbrot set, and infinitely many bounded components;
the latter are simply connected regions with smooth real-algebraic boundary,
and are naturally parameterized by the eigenvalue $\rho \in {\Bbb D}$ of the
unique attracting cycle. Matters become more involved when there are at
least two free critical points. The two-parameter families of normalized
quadratic rational maps and normalized cubic polynomials are often
considered in parallel, as their hyperbolic components admit similar
descriptions: there is a single component of maps with totally disconnected
Julia set and all other components are topological $4$-cells \cite
{milnor2,rees2}. One essential difference is that cubic polynomials with
connected Julia set form a compact set in parameter space; in particular,
every hyperbolic component of maps with two distinct attractors is
precompact. By contrast, while many unbounded hyperbolic components of
quadratic rational maps have been identified \cite{makienko,pilgrim},
bounded components have yet to be exhibited.

Hyperbolic components may also be discussed in the context of Kleinian
groups and their quotient $3$-manifolds. For finitely generated hyperbolic
groups with connected limit set - those whose quotient has incompressible
boundary - the corresponding hyperbolic component is precompact if and only
the limit set is a Sierpinski carpet: the complement of a countable dense
union of Jordan domains with disjoint closures whose diameters tend to zero.
Guided by Sullivan's dictionary between these subjects, McMullen conjectured
that hyperbolic rational maps with Sierpinski carpet Julia set lie in
bounded hyperbolic components \cite{curt}. Pilgrim has suggested more
precisely that a hyperbolic component is bounded when the Julia set is {\em %
almost} a Sierpinski carpet: for example, if every Fatou component is a
Jordan domain and no two Fatou components have closures which intersect in
more than one periodic point. Here we establish precompactness for
hyperbolic components of quadratic rational maps with two attracting cycles,
provided that neither attractor is a fixed point. While it is known in this
case that every Fatou component is a Jordan domain \cite{pilgrim2}, our
largely algebraic arguments do not exploit the topology of the Julia set.

We begin in Section \ref{indsec} with a review of the theory of the
holomorphic index. The index formula 
\[
\frac{1}{1-\alpha }+\frac{1}{1-\beta }+\frac{1}{1-\gamma }=1 
\]
relating the eigenvalues of the three fixed points is fundamental to
Milnor's description \cite{milnor1} of the moduli space of quadratic
rational maps. We survey this work in Section \ref{paramspace} and show in
particular that a sequence of maps is bounded in moduli space if and only if
there is an upper bound on the eigenvalues of the fixed points. Moduli space
is readily parametrized through the choice of a normal form. For certain
purposes it is convenient to work with the family 
\[
f_{\alpha ,\beta }(z)=z\frac{(1-\alpha )z+\alpha (1-\beta )}{\beta (1-\alpha
)z+(1-\beta )} 
\]
of maps fixing $0$,$\infty $, and $1$ with eigenvalues $\alpha $, $\beta $,
and $\gamma =\frac{2-(\alpha +\beta )}{1-\alpha \beta }$; in other settings
it is more useful to work with the family 
\[
F_{\gamma ,\delta }(z)=\frac{\gamma z}{z^{2}+\delta z+1} 
\]
of maps with critical points $\pm 1$ and a fixed point at $0$ with
eigenvalue $\gamma $.

In Section \ref{estsec} we study the limiting dynamics of unbounded
sequences in moduli space. Milnor showed that such sequences accumulate at a
restricted set of points on a natural infinity locus \cite{milnor1},
provided that there are cycles with the same period $n>1$ and uniformly
bounded eigenvalues. We sharpen this and related observations in order to
show that suitably normalized iterates take limits in the family 
\[
G_{T}(z)=z+T+\frac{1}{z} 
\]
as anticipated by considerations in the thesis of Stimson \cite{stim}.
Cycles with bounded eigenvalue tend in the limit to cycles of $G_{T}$ or to
points in the backward orbit of the parabolic fixed point at $\infty $; in
the latter case this backward orbit contains a critical point. In
particular, if the maps in the sequence lie in a hyperbolic component where
there are two nonfixed attractors then $G_{T}$ must have either two
nonrepelling cycles, one nonrepelling cycle and one preperiodic critical
point, or two preperiodic critical points, in addition to the parabolic
fixed point at $\infty $. As discussed in Section \ref{precom}, this
violation of the Fatou-Shishikura bound on the number of nonrepelling cycles
yields the desired contradiction.

Section \ref{conclusion} gives an intersection-theoretic reinterpretation
based on Milnor's observation that ${\rm Per}_n(\rho)$, the locus of
conjugacy classes of maps with an $n$-cycle of eigenvalue $\rho$, is an
algebraic curve whose degree depends only on $n$. The explicit formulas in 
\cite{milnor1} yield a short independent proof of boundedness in the special
case of maps with one attracting cycle of period 2 and another of period 3.
These considerations suggest a combinatorial expression for the intersection
cycle at infinity of a pair of such curves.

\subsection*{Acknowledgments.}

The author wishes to thank Jack Milnor and Kevin Pilgrim for sharing their
insight. Further thanks are due to Misha Yampolsky for assistance with the
graphics which were graciously provided by Jack Milnor. This research has
been generously supported by the Institute for Mathematical Sciences at
Stony Brook.

\section{Local Invariants}

\label{indsec}

Let $g$ be analytic on $U\subseteq {\Bbb C}$ and $\zeta \in U$ with $g(\zeta
)=\zeta $. Assuming that $g$ is not the identity, the {\em topological
multiplicity} is defined as the positive integer 
\[
{\rm mult}_{g}(\zeta )=\frac{1}{2\pi i}\int_{\Gamma }\frac{1-g^{\prime }(z)}{%
z-g(z)}\,dz 
\]
where $\Gamma $ is any sufficiently small positively oriented rectifiable
Jordan curve enclosing $\zeta $; the {\em holomorphic index} is similarly
defined as the complex number 
\[
{\rm ind}_{g}(\zeta )=\frac{1}{2\pi i}\int_{\Gamma }\frac{1}{z-g(z)}\,dz. 
\]
One easily checks that these quantities are invariant under holomorphic
change of coordinates and can thereby be sensibly defined for $\zeta =\infty 
$; moreover, ${\rm mult}_{g}(\zeta )=1$ if and only if the {\em eigenvalue} $%
\rho =g^{\prime }(\zeta )$ differs from $1$, and then 
\begin{equation}
{\rm ind}_{g}(\zeta )=\frac{1}{1-\rho }.  \label{indeig}
\end{equation}
Furthermore, if $|\rho |\neq 1$ or $\rho =1$ then ${\rm mult}_{g^{n}}(\zeta
)={\rm mult}_{g}(\zeta )$ for every $n\geq 1$.

It follows from Cauchy's Integral Formula that 
\begin{eqnarray*}
\sum_{\zeta =g(\zeta )\in V}{\rm mult}_{g}(\zeta ) &=&\frac{1}{2\pi i}%
\int_{\partial V}\frac{1-g^{\prime }(z)}{z-g(z)}\,dz \\
\sum_{\zeta =g(\zeta )\in V}{\rm ind}_{g}(\zeta ) &=&\frac{1}{2\pi i}%
\int_{\partial V}\frac{1}{z-g(z)}\,dz
\end{eqnarray*}
for open $V$ with $\overline{V}\subseteq U\subseteq {\Bbb C}$ and with
rectifiable boundary containing no fixed points. These sums evidently depend
continuously on $g$. For rational maps $g:\widehat{{\Bbb C}}\rightarrow 
\widehat{{\Bbb C}}$ of degree $d$, one sees from the Residue Theorem that 
\begin{equation}
\sum_{\zeta =g(\zeta )\in \hat{{\bf C}}}{\rm mult}_{g}(\zeta )=d+1;
\label{topmult}
\end{equation}
the {\em Holomorphic Index Formula} 
\begin{equation}
\sum_{\zeta =g(\zeta )\in \hat{{\bf C}}}{\rm ind}_{g}(\zeta )=1
\label{holindex}
\end{equation}
follows similarly. We denote ${\rm Fix}(g)$ the unordered $(d+1)$-tuple of
fixed points listed with multiplicity. In general, we denote such
collections of possibly identical points as $[x_{1},\ldots ,x_{n}]$. We
similarly write ${\rm Crit}(g)$ for the unordered $(2d-2)$-tuple of critical
points; note that there are at least two distinct critical points when $%
d\geq 2$.

A fixed point $\zeta $ of an analytic map $g$ is said to be {\em attracting,
indifferent}, or {\em repelling} according as the eigenvalue $\rho $ is less
than, equal to, or greater than $1$. If $\rho =e^{2\pi ip/q}$ where $(p,q)=1$
and $g^{q}$ is not the identity, then $\zeta $ is {\em parabolic}. A simple
calculation then shows that ${\rm mult}_{g^{q}}(\zeta )=\ell q+1$ for some
positive integer $\ell $; we refer to $\ell $ as the {\em degeneracy}, and
say that $\zeta $ is a {\em degenerate parabolic} fixed point when $\ell
\geq 2$. In view of (\ref{indeig}), if ${\rm mult}_{g}(\zeta )=1$ then $%
\zeta $ is attracting, indifferent, or repelling according as the real part
of $\,{\rm ind}_{g}(\zeta )$ is greater than, equal to, or less than $\frac{1%
}{2}$. Following \cite{algdyn} we say that a parabolic fixed point $\zeta $
with eigenvalue $e^{2\pi ip/q}$ is 
\[
\begin{array}{lll}
{\it parabolic-attracting}\;\; & \mbox{ when }\;\; & \Re \,{\rm ind}%
_{g^{q}}(\zeta )>\frac{\ell q+1}{2}, \\ 
{\it parabolic-indifferent}\;\; & \mbox{ when }\;\; & \Re \,{\rm ind}%
_{g^{q}}(\zeta )=\frac{\ell q+1}{2}, \\ 
{\it parabolic-repelling}\;\; & \mbox{ when }\;\; & \Re \,{\rm ind}%
_{g^{q}}(\zeta )<\frac{\ell q+1}{2}.
\end{array}
\]

More generally, we say that $\zeta $ is {\em periodic} under $g$ when $%
g^{n}(\zeta )=\zeta $ for some $n\geq 1$, the least such $n$ being referred
to as the {\em period}. The multiplicity, index, and eigenvalue of the cycle 
$\langle \zeta \rangle =\{\zeta ,\ldots ,g^{n-1}(\zeta )\}$ are the
corresponding invariants of $\zeta $ as a fixed point of $g^{n}$. It follows
from the definition of multiplicity that a generic perturbation of $g$
splits an $n$-cycle with eigenvalue $\rho =e^{2\pi ip/q}$ and degeneracy $%
\ell $ into an $n$-cycle with eigenvalue close to $\rho $ and an $\ell $%
-tuple of $nq$-cycles with eigenvalues close to $1$. Continuity of the local
index sum implies:

\begin{lemma}
\label{paratt} Let $g$ be analytic on $U$ with a parabolic $n$-cycle $%
\langle \zeta \rangle $ of eigenvalue $e^{2\pi ip/q}$. Further let $g_{k}$
be analytic with $g_{k}\rightarrow g$ locally uniformly on $U$, and with $n$%
-cycles $\langle \zeta _{k}^{[0]}\rangle $ and $nq$-cycles $\langle \zeta
_{k}^{[1]}\rangle ,\ldots ,\langle \zeta _{k}^{[\ell ]}\rangle $ converging
to $\langle \zeta \rangle $. If all $\langle \zeta _{k}^{[j]}\rangle $ are
attracting for $k$ sufficiently large then $\langle \zeta \rangle $ is
parabolic-attracting or parabolic-indifferent.
\end{lemma}

\smallskip Assume now that $g$ is rational of degree $d$. The {\em basin} of
an attracting cycle $\langle \zeta \rangle $ is the open set consisting of
all points $z\in \widehat{{\Bbb C}}$ with $g^{n}(z)\rightarrow \langle \zeta
\rangle $. We refer to the connected component containing $\xi \in \langle
\zeta \rangle $ as the {\em immediate basin} of $\xi $. The basin of a
parabolic cycle is similarly defined as the open set of all $z\in \widehat{%
{\Bbb C}}$ with $\langle \zeta \rangle \not{\ni}g^{n}(z)\rightarrow \langle
\zeta \rangle $, the $\ell q$ components adjoining $\xi $ forming the
immediate basin of $\xi $. In both cases, the immediate basin of $\langle
\zeta \rangle $ is taken to be the union of the immediate basins of the
points in the cycle. Fatou established the fundamental fact that each cycle
of components of the immediate basin of an attracting or parabolic cycle
always contains at least one critical value with infinite forward orbit \cite
{milnor0}. In particular, counting degeneracy there are at most $2d-2$
attracting and parabolic cycles. Shishikura extended this bound to the total
count of nonrepelling cycles \cite{shish}, and the author proved a refined
inequality where the contribution of each parabolic-attracting and
parabolic-indifferent cycle is augmented by one \cite{algdyn}; consideration
of the return maps on Ecalle cylinders shows in fact that there are at least 
$\ell +1$ critical values with infinite forward orbit in the immediate basin
of a parabolic-attracting or parabolic-indifferent cycle of degeneracy $\ell 
$.

Consider the family 
\[
G_{T}(z)=z+T+\frac{1}{z} 
\]
of quadratic rational maps with critical points $\pm 1$ and a degenerate
fixed point at $\infty $ with eigenvalue $1$ and holomorphic index $1-\frac{1%
}{T^{2}}$; by convention, $G_{\infty }\equiv \infty $. The Fatou-Shishikura
Inequality has the following consequences in this special case:

\begin{lemma}
\label{fsineq} Let $G=G_{T}$ where $T\in {\Bbb C}$.

\begin{itemize}
\item  If $T=0$ then $\infty $ is a degenerate parabolic fixed point.
Neither critical point is preperiodic and all other cycles are repelling.

\item  If $T\neq 0$ and $\langle \zeta \rangle $ is attracting or
indifferent, then neither critical point is preperiodic and all other cycles
are repelling; if parabolic, then $\langle \zeta \rangle $ is nondegenerate
parabolic-repelling.

\item  If $T\neq 0$ and either critical point is preperiodic, then the other
critical point has infinite forward orbit and all other cycles are repelling.
\end{itemize}
\end{lemma}

\begin{figure}[tbp]
\centerline{\psfig{figure=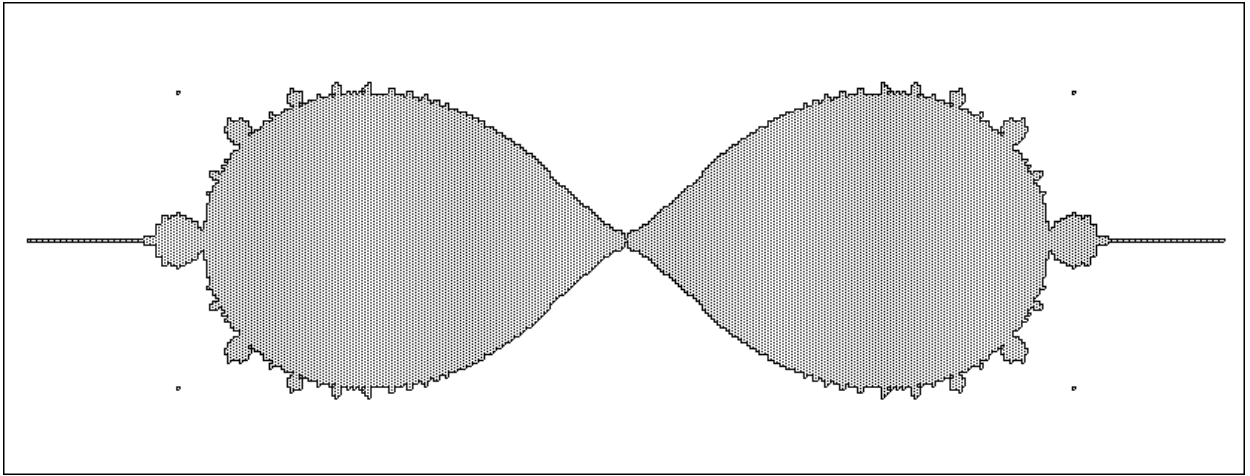,width=\hsize}}
\caption{Bifurcation locus for the family $G_\kappa(z)=z+\kappa+ \frac{1}{z}$%
.}
\label{gkfig}
\end{figure}

\section{Normal Forms}

\label{paramspace}

We naturally identify the space of all quadratic rational maps 
\[
{\cal RAT}_{2}=\left\{ g(z)=\frac{A_{2}z^{2}+A_{1}z+A_{0}}{%
B_{2}z^{2}+B_{1}z+B_{0}}:\;\deg g=2\right\} 
\]
with the open subvariety of projective space ${\Bbb P}^{5}$ where the
resultant 
\[
{\rm det}\left( 
\begin{array}{cccc}
A_{2} & A_{1} & A_{0} & 0 \\ 
0 & A_{2} & A_{1} & A_{0} \\ 
B_{2} & B_{1} & B_{0} & 0 \\ 
0 & B_{2} & B_{1} & B_{0}
\end{array}
\right) 
\]
is nonvanishing. Various technical purposes require that we work in the
spaces 
\begin{eqnarray*}
{\cal RAT}_{2}^{\times } &=&\left\{ (g;\;\chi ^{+},\chi ^{-})\in {\cal RAT}%
_{2}\times \widehat{{\Bbb C}}^{2}:\;{\rm Crit}(g)=[\chi ^{+},\chi
^{-}]\;\right\}  \\
{\cal RAT}_{2}^{\circ } &=&\left\{ (g;\;a,b,c)\in {\cal RAT}_{2}\times 
\widehat{{\Bbb C}}^{3}:\;{\rm Fix}(g)=[a,b,c]\;\right\}  \\
{\cal RAT}_{2}^{\otimes } &=&\left\{ (g;\;\chi ^{+},\chi ^{-};\;a,b,c)\in 
{\cal RAT}_{2}\times \widehat{{\Bbb C}}^{5}:
\begin{array}{cc}
{\rm Crit}(g)=[\chi ^{+},\chi ^{-}] &  \\ 
{\rm Fix}(g)=[a,b,c] & 
\end{array}
\right\} 
\end{eqnarray*}
where the critical points, fixed points, or both have been marked. The
quotients under the conjugation action of the M\"{o}bius group are the {\em %
moduli spaces} 
\begin{eqnarray*}
{\bf Rat}_{2} &=&{\cal RAT}_{2}/{\rm PSL}_{2}{\Bbb C} \\
{\bf Rat}_{2}^{\times } &=&{\cal RAT}_{2}^{\times }/{\rm PSL}_{2}{\Bbb C} \\
{\bf Rat}_{2}^{\circ } &=&{\cal RAT}_{2}^{\circ }/{\rm PSL}_{2}{\Bbb C} \\
{\bf Rat}_{2}^{\otimes } &=&{\cal RAT}_{2}^{\otimes }/{\rm PSL}_{2}{\Bbb C}
\end{eqnarray*}
all varieties of complex dimension 2.

Writing $\alpha,\beta,\gamma$ for the eigenvalues of the fixed points $a,b,c$
we see from (\ref{holindex}) that 
\[
\frac{1}{1-\alpha}+\frac{1}{1-\beta}+\frac{1}{1-\gamma}=1 
\]
so long as $\alpha,\beta,\gamma\neq 1$, and 
\begin{equation}  \label{abc}
\alpha\beta\gamma-(\alpha+\beta+\gamma)+2=0
\end{equation}
always; in particular, 
\[
\gamma=\frac{2-(\alpha+\beta)}{1-\alpha\beta}. 
\]

Let $[(g;\;\chi ^{+},\chi ^{-}\;;a,b,c)]]$ be a class in ${\bf Rat}%
_{2}^{\otimes }$. Provided that $\chi ^{+}\neq c\neq \chi ^{-}$, there is a
unique representative of the form $(F;\;+1,-1;\;a,b,0)$ where 
\begin{equation}
F(z)=F_{\gamma ,\delta }(z)=\frac{\gamma z}{z^{2}+\delta z+1}  \label{Fnf}
\end{equation}
for some $\gamma ,\delta \in {\Bbb C}$ with $\gamma \neq 0$; moreover, every
class in ${\bf Rat}_{2}^{\times }$ has a representative of this form. As 
\[
F_{\gamma ,\delta }^{\prime }(z)=\frac{\gamma (1-z^{2})}{z^{2}+\delta z+1} 
\]
it follows that 
\begin{eqnarray*}
\alpha =\frac{1-a^{2}}{\gamma } &=&\frac{\delta a+2}{\gamma }-1 \\
\beta =\frac{1-b^{2}}{\gamma } &=&\frac{\delta b+2}{\gamma }-1
\end{eqnarray*}
with 
\[
\{a,b\}=\left\{ \frac{-\delta \pm \sqrt{\delta ^{2}-4(1-\gamma )}}{2}%
\right\} . 
\]

Alternatively, provided that $a\neq b\neq c\neq a$ there is a unique
representative of the form 
\begin{equation}
f_{\alpha ,\beta }(z)=z\frac{(1-\alpha )z+\alpha (1-\beta )}{\beta (1-\alpha
)z+(1-\beta )}  \label{fnf}
\end{equation}
for some $\alpha ,\beta \in {\Bbb C}$ with $\alpha ,\beta ,\alpha \beta \neq
1$. Writing 
\begin{eqnarray}
f_{\alpha ,\beta }(z) &=&z\frac{(1-\alpha )(z-1)+\epsilon }{\beta (1-\alpha
)(z-1)+\epsilon }  \label{fab1} \\
&=&\frac{z}{\beta }\left[ \frac{z-\nu }{z-\mu }\right] =\frac{z}{\beta }%
\left[ 1+\frac{\mu \epsilon }{z-\mu }\right] ,  \label{fab2}
\end{eqnarray}
where 
\begin{equation}
\epsilon =1-\alpha \beta =\frac{(1-\alpha )(1-\beta )}{(\gamma -1)}
\label{epsdef}
\end{equation}
and 
\begin{equation}
\mu =\frac{\beta -1}{\beta -\alpha \beta }=1-\frac{\epsilon }{\beta
(1-\alpha )}  \label{mudef}
\end{equation}
\begin{equation}
\nu =\frac{\alpha \beta -\alpha }{1-\alpha }=1-\frac{\epsilon }{1-\alpha }\,,
\label{nudef}
\end{equation}
we see that $f_{\alpha ,\beta }(\mu )=\infty $ and $f_{\alpha ,\beta }(\nu
)=0$. Calculating the derivative 
\begin{eqnarray}
f_{\alpha ,\beta }^{\prime }(z) &=&\frac{\beta (1-\alpha
)^{2}(z-1)^{2}+(1+\beta )(1-\alpha )(z-1)\epsilon +(2-\alpha -\beta
)\epsilon }{[\beta (1-\alpha )(z-1)+\epsilon ]^{2}}  \label{dfab1} \\
&=&\frac{1}{\beta }\left[ 1+\frac{\mu \nu -\mu ^{2}}{(z-\mu )^{2}}\right]
\;=\;\frac{1}{\beta }\left[ 1-\frac{\mu ^{2}\epsilon }{(z-\mu )^{2}}\right]
\,,  \label{dfab2}
\end{eqnarray}
we find that 
\[
\mu =\frac{\chi ^{+}+\chi ^{-}}{2} 
\]
\[
\epsilon =\left( \frac{\chi ^{+}-\chi ^{-}}{\chi ^{+}+\chi ^{-}}\right) ^{2} 
\]
whence 
\[
\chi ^{\pm }=\mu (1\pm \sqrt{\epsilon }) 
\]
for the appropriate choice of $\sqrt{\epsilon }$.

Assuming both restrictions on the marked points, there is a unique
M\"{o}bius transformation 
\begin{equation}
\phi (z)=\frac{bz-ab}{az-ab}=\frac{(\mu ^{2}\epsilon -\mu ^{2}+\mu )z+\mu 
\sqrt{\epsilon }}{(1-\mu )z+\mu \sqrt{\epsilon }}  \label{phiz}
\end{equation}
sending $+1,-1,a,b,0$ to $\chi ^{+},\chi ^{-},0,\infty ,1$. Clearly, 
\[
F_{\gamma ,\delta }=\phi ^{-1}\circ f_{\alpha ,\beta }\circ \phi 
\]
where 
\[
(\alpha ,\beta )=\left( \frac{4\chi ^{+}\chi ^{-}-2\chi ^{+}\chi ^{-}(\chi
^{+}+\chi ^{-})}{(\chi ^{+}+\chi ^{-})^{2}-2\chi ^{+}\chi ^{-}(\chi
^{+}+\chi ^{-})}\,,\,\frac{2(\chi ^{+}+\chi ^{-})-4\chi ^{+}\chi ^{-}}{%
2(\chi ^{+}+\chi ^{-})-(\chi ^{+}+\chi ^{-})^{2}}\right) 
\]
and 
\[
(\gamma ,\delta )=(1-ab\,,\,-a-b). 
\]

Recall that the elementary symmetric functions 
\begin{eqnarray*}
X(\alpha ,\beta ,\gamma ) &=&\alpha +\beta +\gamma \\
Y(\alpha ,\beta ,\gamma ) &=&\alpha \beta +\alpha \gamma +\beta \gamma \\
Z(\alpha ,\beta ,\gamma ) &=&\alpha \beta \gamma
\end{eqnarray*}
together determine $[\alpha ,\beta ,\gamma ]$. It follows from (\ref{abc})
that 
\[
{\bf Rat}_{2}^{\circ }\ni [f;\;a,b,c]\leadsto (X(\alpha ,\beta ,\gamma
),Y(\alpha ,\beta ,\gamma ),Z(\alpha ,\beta ,\gamma ))\in {\Bbb C}^{3} 
\]
descends to a map ${\bf Rat}_{2}\rightarrow {\Bbb C}^{3}$ with image in the
hyperplane 
\[
\{(X,Y,Z)\in {\Bbb C}^{3}:\;Z=X-2\}, 
\]
and we obtain 
\[
j:{\bf Rat}_{2}\rightarrow {\Bbb C}^{2} 
\]
on composing with the projection ${\Bbb C}^{3}\ni (X,Y,Z)\leadsto (X,Y)\in 
{\Bbb C}^{2}$. Consideration of the normal forms (\ref{Fnf}) and (\ref{fnf})
shows that an unordered triple $[\alpha ,\beta ,\gamma ]$ satisfying (\ref
{abc}) determines a unique class in ${\bf Rat}_{2}$, and thus $j$ is an
isomorphism. As $\{(X(\alpha_{k},\beta _{k},\gamma _{k}),Y(\alpha _{k},\beta
_{k},\gamma _{k})\}$ and $\{\alpha _{k},\beta _{k},\gamma _{k}\}$ are 
simultaneously bounded or unbounded, we recover Milnor's observation \cite
{milnor1}:

\begin{lemma}
\label{unbdd} Let $g_{k}$ be quadratic rational maps with eigenvalues $%
\alpha _{k},\beta _{k},\gamma _{k}$ at the fixed points $a,b,c$. Then $%
[g_{k}]$ is bounded in ${\bf Rat}_{2}$ if and only if $\{\alpha _{k},\beta
_{k},\gamma _{k}\}$ is bounded in ${\Bbb C}$.
\end{lemma}

\section{Limit Dynamics}

\label{estsec}

Our first goal is the following:

\begin{proposition}
\label{no01} Let $g_{k}$ be quadratic rational maps with eigenvalues $\alpha
_{k},\beta _{k},\gamma _{k}$ at the fixed points $a,b,c$, where $\alpha _{k}$
and $\beta _{k}$ converge in $\widehat{{\Bbb C}}$ and $\gamma
_{k}\rightarrow \infty $. Assume that there are cycles $\langle z_{k}\rangle 
$ with the same period $n>1$ and uniformly bounded eigenvalues. Then 
\[
\alpha _{k}=\omega +O(\sqrt{\epsilon _{k}}) \text{
\qquad and \qquad }\beta _{k}=\bar{\omega}+O(\sqrt{\epsilon_{k}})\]
as $k\rightarrow \infty $, where $\omega \neq 1$ is a $q$-th root of unity
for some $q\leq n$ and 
\[\epsilon_k=1-\alpha_k\beta_k=O\left(\frac{1}{\gamma_k}\right).\]
\end{proposition}

\medskip The proof requires several preliminary lemmas and the elimination
of various special cases. Let $\alpha _{k},\beta _{k},\gamma _{k}\in {\Bbb C}
$ satisfying $(\ref{abc})$, and suppose that $\gamma _{k}\rightarrow \infty $%
. Inspection of $(\ref{epsdef})$, $(\ref{mudef})$, and $(\ref{nudef})$ shows
that

\begin{itemize}
\item  $\alpha _{k}\rightarrow \infty $ if and only if $\beta
_{k}\rightarrow 0$, and vice-versa;

\item  $\alpha _{k}\rightarrow 1$ if and only if $\beta _{k}\rightarrow 1$,
and vice-versa;

\item  $\epsilon _{k}=o(\alpha _{k}-1)$ if $\beta _{k}$ is bounded, and $%
\epsilon _{k}=o(\beta _{k}-1)$ if $\alpha _{k}$ is bounded;

\item  $\epsilon _{k}=O(\gamma_k^{-1})$ if both $\alpha _{k}$ and $\beta _{k}$
are bounded;

\item  $\mu _{k}\rightarrow 1$ if $\alpha _{k}\neq 1$ is bounded, and $\nu
_{k}\rightarrow 1$ if $\beta _{k}\neq 1$ is bounded;

\item  Both $\mu _{k}$ and $\nu _{k}$ are $1+O(\epsilon _{k})$ if $\alpha
_{k}$, hence also $\beta _{k}$, is bounded away from $\{0,1,\infty \}$.
\end{itemize}

Recall from (\ref{phiz}) that the choice of $\sqrt{\epsilon _{k}}$,
corresponding to a marking of the critical points, specifies a M\"{o}bius
transformation $\phi _{k}$ which conjugates $f_{\alpha _{k},\beta _{k}}$ to
some $F_{\gamma _{k},\delta _{k}}$. It follows from these observations that 
\begin{equation}
\phi _{k}(z)=1+z\sqrt{\epsilon _{k}}+o(\sqrt{\epsilon _{k}})  \label{magnify}
\end{equation}
on compact sets in ${\Bbb C}$, provided that $\alpha _{k}$ and $\beta _{k}$
are bounded away from $\{0,1,\infty \}$.

\bigskip Let $f_{k}=f_{\alpha _{k},\beta _{k}}$ where $\gamma
_{k}\rightarrow \infty $ and $\alpha _{k}\rightarrow \alpha _{\infty }\in 
{\Bbb C}^{*}$. By (\ref{fab1}) and (\ref{fab2}), 
\[
\frac{f_{k}(z_{k})}{z_{k}}=1+o(1)\;\;\mbox{ if }z_{k}-1=o(\epsilon _{k}) 
\]
and 
\[
f_{k}(z)\rightarrow \alpha _{\infty }z\;\;\mbox{ locally uniformly on }%
\widehat{{\Bbb C}}{\bf -\{}0,1,\infty \} 
\]
Assuming further that $\alpha _{\infty }\neq 1$, we have 
\[
\frac{f_{k}(z_{k})}{z_{k}}=\alpha _{k}\left[ 1+\frac{\epsilon _{k}}{z_{k}-1}%
+o\left( \frac{\epsilon _{k}}{z_{k}-1}\right) +O(\epsilon _{k})\right] 
\]
when $\epsilon _{k}=o(z_{k}-1)$, and thus 
\begin{equation}
\frac{f_{k}(z_{k})}{z_{k}}=\left\{ 
\begin{array}{ll}
\alpha _{k}+o(1) & \mbox{ if }\epsilon _{k}=o(z_{k}-1) \\ 
\alpha _{k}(1+\frac{1}{\tau }\sqrt{\epsilon _{k}})+o(\sqrt{\epsilon _{k}}) & %
\mbox{ if }z_{k}=1+\tau \sqrt{\epsilon _{k}}+o(\sqrt{\epsilon _{k}})%
\mbox{ for
}\tau \in {\Bbb C}^{*} \\ 
\alpha _{k}+o(\sqrt{\epsilon _{k}}) & \mbox{ if }\sqrt{\epsilon _{k}}%
=o(z_{k}-1) \\ 
\alpha _{k}+O(\epsilon _{k}) & \mbox{ if }z_{k}\mbox{ is bounded away from }%
1;
\end{array}
\right.  \label{foverz}
\end{equation}
moreover, 
\begin{equation}
\frac{f_{k}(z_{k})-1}{z_{k}-1}=\frac{(z_{k}-\beta _{k})(1-\alpha
_{k})+\epsilon _{k}}{\beta _{k}(1-\alpha _{k})(z_{k}-1)+\epsilon _{k}}%
\rightarrow \infty  \label{f1z1}
\end{equation}
whenever $z_{k}\rightarrow 1$.

\medskip Observe that if $z_{k}=1+O(\sqrt{\epsilon _{k}})$ and $%
f_{k}^{n}(z_{k})=1+O(\sqrt{\epsilon _{k}})$ then 
\[
1+O(\sqrt{\epsilon _{k}})=\frac{f_{k}^{n}(z_{k})}{z_{k}}=\prod_{j=0}^{n-1}%
\frac{f_{k}^{j+1}(z_{k})}{f_{k}^{j}(z_{k})}=\alpha _{k}^{n}+O(\sqrt{\epsilon
_{k}}), 
\]
unless $f_{k}^{j}(z_{k_{\ell }})=1+o(\sqrt{\epsilon _{k}})$ for some $0\leq
j<n$ and $k_{\ell }\rightarrow \infty $. Applying (\ref{magnify}) we deduce:

\begin{lemma}
\label{alphan} Let $F_{k}=F_{\gamma _{k},\delta _{k}}$ where $\gamma
_{k}\rightarrow \infty $ and $\alpha _{k}\rightarrow \alpha _{\infty }\notin
\{0,1,\infty \}$, hence $\beta _{k}\rightarrow \beta _{\infty }=\alpha
_{\infty }^{-1}.$ Suppose that $z_{k}\in \widehat{{\Bbb C}}$ with $%
F_{k}^{j}(z_{k})\rightarrow \zeta ^{(j)}\in \widehat{{\Bbb C}}$ for $0\leq
j\leq n$, where $\zeta ^{(0)}\not{\in}\{0,\infty \}$ and $\zeta ^{(n)}\neq
\infty $. If $\zeta ^{(j)}\neq 0$ for $0<j<n$ then 
\[
\alpha _{k}=\alpha _{\infty }+O(\sqrt{\epsilon _{k}})\text{ \qquad and
\qquad }\beta _{k}=\beta _{\infty }+O(\sqrt{\epsilon _{k}}) 
\]
and $\alpha _{\infty }^{n}=1=\beta _{\infty }^{n}$; moreover, if $\zeta
^{(j)}=\infty $ for $0<j<n$ then $\alpha _{\infty }$, hence also $\beta
_{\infty },$ is a primitive $n$-th root of unity.
\end{lemma}

\smallskip Assume now that $\alpha _{k}\rightarrow \omega $, hence $\beta
_{k}\rightarrow \bar{\omega}$, where $\omega \neq 1$ is a root of unity, and
suppose that $z_{k}=1+o(\sqrt{\epsilon _{k}})$ but $\epsilon _{k}=o(z_{k}-1)$%
. Then 
\[
f_{k}^{n}(z_{k})=\omega ^{n}\left[ 1+\frac{\epsilon _{k}}{z_{k}-1}+o\left( 
\frac{\epsilon _{k}}{z_{k}-1}\right) \right] 
\]
for $n\geq 1$ by (\ref{foverz}) and induction, and thus (\ref{magnify})
implies:

\begin{lemma}
\label{scaledown} Let $F_{k}=F_{\gamma _{k},\delta _{k}}$ where $\gamma
_{k}\rightarrow \infty $. Assume that $\alpha _{k}\rightarrow \omega $,
hence $\beta _{k}\rightarrow \bar{\omega}$, where $\omega \neq 1$ is a root
of unity, and let $z_{k}\in \widehat{{\Bbb C}}$ with $z_{k}\rightarrow 0$.
If $F_{k}^{n}(z_{k})$ is bounded for some $n>1$ then $z_{k}=O(\sqrt{\epsilon
_{k}})$.
\end{lemma}

Suppose now that $f_{k}^{n}(z_{k})=z_{k}$ where $z_{k}\in \widehat{{\Bbb C}}%
-\{0,1,\infty \}$ and $n>1$. As 
\[
1=\frac{f_{k}^{n}(z_{k})-1}{z_{k}-1}=\prod_{j=0}^{n-1}\frac{%
f_{k}^{j+1}(z_{k})-1}{f_{k}^{j}(z_{k})-1}, 
\]
it follows from (\ref{f1z1}) that $\zeta _{k}$ is bounded away from $1$ for
some $\zeta _{k}\in \langle z_{k}\rangle $. Similarly, 
\[
\alpha _{k}^{n}=\left\{ 
\begin{array}{ll}
1+o(1) & \mbox{ if }\min_{\zeta \in \langle z_{k}\rangle }\left| \frac{\zeta
-1}{\epsilon _{k}}\right| \rightarrow \infty \\ 
&  \\ 
1+O(\sqrt{\epsilon _{k}}) & \mbox{ if }\min_{\zeta \in \langle z_{k}\rangle
}\left| \frac{\zeta -1}{\sqrt{\epsilon _{k}}}\right| 
\mbox{ is bounded away
from $0$} \\ 
&  \\ 
1+o(\sqrt{\epsilon _{k}}) & \mbox{ if }\min_{\zeta \in \langle z_{k}\rangle
}\left| \frac{\zeta -1}{\sqrt{\epsilon _{k}}}\right| \rightarrow \infty \\ 
&  \\ 
1+O(\epsilon _{k}) & \mbox{ if }\langle z_{k}\rangle 
\mbox{ is bounded away
from }1.
\end{array}
\right. 
\]
by (\ref{foverz}). Combining these observations with Lemma \ref{alphan}, we
obtain:

\begin{lemma}
\label{calz}Let $F_{k}=F_{\gamma _{k},\delta _{k}}$ where $\gamma
_{k}\rightarrow \infty $. Assume that $\alpha _{k}\rightarrow \omega $,
hence $\beta _{k}\rightarrow \bar{\omega}$, where $\omega \neq 1$ is a root
of unity, and let $\langle z_{k}\rangle $ be cycles of period $n>1$. If $%
\langle z_{k}\rangle \rightarrow {\cal Z}\subset \widehat{{\Bbb C}}$ then $%
\infty \in {\cal Z}$. Moreover:

\begin{itemize}
\item  If ${\cal Z}\neq \{0,\infty \}$ then $\alpha _{k}=\omega +O(\sqrt{%
\epsilon _{k}})$ and $\beta _{k}=\bar{\omega}+O(\sqrt{\epsilon _{k}}).$

\item  If ${\cal Z}=\{\infty \}$ then $\alpha _{k}=\omega +o(\sqrt{\epsilon
_{k}})$ and $\beta _{k}=\bar{\omega}+o(\sqrt{\epsilon _{k}}).$
\end{itemize}
\end{lemma}

\medskip Assume now that $\alpha _{k}$, hence also $\beta _{k}$, is bounded
away from $\{0,1,\infty \}$. It follows from (\ref{dfab1}) that 
\[
f_{k}^{\prime }(z_{k})=\frac{(2-\alpha _{k}-\beta _{k})\epsilon
_{k}+o(\epsilon _{k})}{o(\epsilon _{k})}\rightarrow \infty 
\]
if $z_{k}-1=o(\sqrt{\epsilon _{k}})$. On the other hand, if $\sqrt{\epsilon
_{k}}=O(z_{k}-1)$ then (\ref{dfab2}) implies 
\[
f_{k}^{\prime }(z_{k})=\frac{1}{\beta _{k}}\left[ 1-\frac{\epsilon _{k}}{%
(z_{k}-1)^{2}}+o\left( \frac{\epsilon _{k}}{(z_{k}-1)^{2}}\right) \right] 
\]
whence 
\[
f_{k}^{\prime }(z_{k})=\left\{ 
\begin{array}{ll}
\alpha _{k}\left( 1-\frac{1}{\tau ^{2}}\right) +o(1) & \mbox{ if }%
z_{k}=1+\tau \sqrt{\epsilon _{k}}+o(\sqrt{\epsilon _{k}})\mbox{ for }\tau
\in {\Bbb C}^{*}, \\ 
\alpha _{k}+o(1) & \mbox{ if }\sqrt{\epsilon _{k}}=o(z_{k}-1).
\end{array}
\right. 
\]
In particular:

\begin{lemma}
\label{pqbddeig} Let $F_{k}=F_{\gamma _{k},\delta _{k}}$ where $\gamma
_{k}\rightarrow \infty $. Assume that $\alpha _{k}\rightarrow \omega $,
hence $\beta _{k}\rightarrow \bar{\omega}$, where $\omega \neq 1$ is a
primitive $q$-th root of unity, and let $\langle z_{k}\rangle \rightarrow 
{\cal Z}$ be cycles of period $n>1$ and eigenvalues $\rho _{k}$. If $0\not%
{\in}{\cal Z}$ then $q|n$ and $\rho _{k}$ is bounded; moreover, if ${\cal Z}%
=\{\infty \}$ then $\rho _{k}\rightarrow 1$.
\end{lemma}

We similarly deduce:

\begin{lemma}
\label{bddeig}Let $F_{k}=F_{\gamma _{k},\delta _{k}}$ where $\gamma
_{k}\rightarrow \infty $. Assume that $\alpha _{k}$, hence also $\beta _{k}$%
, is bounded away from $\{0,1,\infty \}$, and let $\langle z_{k}\rangle
\rightarrow {\cal Z}$ be cycles of period $n>1$ and eigenvalues $\rho _{k}$.
If $0\in {\cal Z}$ and $\rho _{k}$ is bounded then $+1\in {\cal Z}$ or $%
-1\in {\cal Z}$.
\end{lemma}

\bigskip \noindent {\bf Proof of Proposition \ref{no01}:} Assume that $%
\alpha _{k}\rightarrow \alpha _{\infty }$. If $\alpha _{\infty }\not{\in}\
\{0,1,\infty \}$ then also $\beta _{\infty }\not{\in}\{0,1,\infty \}$. In
particular, we may represent each class 
\[
\lbrack (g_{k};\;a_{k},b_{k},c_{k})]]\in {\bf Rat}_{2}^{\circ } 
\]
by a map $F_{k}=F_{\gamma _{k},\delta _{k}}$. Recall that we are given $n$%
-cycles $\langle z_{k}\rangle $ with uniformly bounded eigenvalues $\rho
_{k} $. Passing to a subsequence if necessary, we may assume that $\langle
z_{k}\rangle \rightarrow {\cal Z}\subseteq \hat{{\bf C}}$. In view of Lemma 
\ref{bddeig}, either ${\cal Z}=\{\infty \}$ or ${\cal Z}\cap {\Bbb C}%
^{*}\neq \emptyset $, and the conclusion
follows by Lemma \ref{calz}.

Suppose next that $\alpha _{\infty }$, hence also $\beta _{\infty }$, is in $%
\{0,\infty \}$. Permuting the fixed points if necessary, we may assume on
passage to a subsequence that $\alpha _{k}\rightarrow \infty $ and $\alpha
_{k}=O(\gamma _{k})$. Following Milnor \cite{milnor1} we work with the
representatives $(\hat{f}_{k};\;0,\infty ,c_{k})$ where 
\[
\hat{f}_{k}(z)=z\frac{z+\alpha _{k}}{\beta _{k}z+1} 
\]
and 
\[
c_{k}=\frac{1-\alpha _{k}}{1-\beta _{k}}=-\alpha _{k}+o(\alpha _{k}). 
\]
Calculating the derivative 
\[
\hat{f}_{k}^{\prime }(z)=\frac{\beta _{k}z^{2}+2z+\alpha _{k}}{(\beta
_{k}z+1)^{2}} 
\]
we see that $\hat{f}_{k}^{\prime }(z)=\alpha _{k}+O(1)$ on the disc $|z|<4$.
In particular, $\hat{f}_{k}$ is univalent on $|z|<4$ with image containing
the disc $|z|<3|c_{k}|$, and both critical values lie outside the latter
region. Consequently, there are univalent inverse branches $A_{k}$ and $%
C_{k} $, fixing $0$ and $c_{k}$, defined on the disc $|z|<3|c_{k}|$. As $%
D_{k}=\{z:|2z-c_{k}|<2|c_{k}|\}$ lies in the image of the disc $|z|<4$, it
follows that $A_{k}^{\prime }(D_{k})=O(\alpha _{k}^{-1})$ on $D_{k}$ and $%
A_{k}(D_{k})\subset D_{k}$. On the other hand, 
\[
C_{k}^{\prime }(z)=O(\gamma _{k}^{-1})=O(\alpha _{k}^{-1}) 
\]
for $|z|<\frac{5}{2}|c_{k}|$ by the compactness of normalized univalent
functions; consequently, $|C_{k}(z)-c_{k}|=O(c_{k}\gamma _{k}^{-1})=O(1)$
for $|z-c_{k}|<\frac{3}{2}|c_{k}|$, and in particular $C_{k}(D_{k})\subset
D_{k}$. We deduce that $J(\hat{f}_{k})\subset \hat{f}_{k}^{-1}(D_{k})$ is a
Cantor set containing all periodic points other than the fixed point at $%
\infty $. Thus, $\langle z_{k}\rangle \subset J(\hat{f}_{k})$ and $\rho
_{k}^{-1}=O(\alpha _{k}^{-n})$ whence $\rho _{k}\rightarrow \infty $.

It remains to treat the case $\alpha _{\infty }=1=\beta _{\infty }$. Now it
is advantageous to choose representatives $(\hat{g}_{k};\;\infty
,b_{k},c_{k})$ where 
\[
\hat{g}_{k}(z)=\frac{(\alpha _{k}\gamma _{k}-1)z^{2}+(\alpha _{k}^{2}\gamma
_{k}-\alpha _{k}^{2})z+\alpha _{k}^{2}}{(\alpha _{k}^{2}\gamma _{k}-\alpha
_{k})z} 
\]
and 
\[
(b_{k},c_{k})=\left( \frac{\alpha _{k}}{\alpha _{k}-1},\frac{\alpha _{k}}{%
1-\alpha _{k}\gamma _{k}}\right) \rightarrow (\infty ,0). 
\]
Notice that $\hat{g}_{k}(z)\rightarrow z+1$ locally uniformly on $\widehat{%
{\Bbb C}}-\{0\}$, and thus $\hat{g}_{k}^{n}(z)\rightarrow z+n$ locally
uniformly on $\widehat{{\Bbb C}}-\{-(n-1),\ldots ,0\}$. As the translation $%
z\leadsto z+1$ has a fixed point of multiplicity $2$ at $\infty $ and no
other periodic points, $b_{k}$ and $\infty $ are the only fixed points of $%
\hat{g}_{k}^{n}$ outside the circle $|z|=n$. We may therefore assume without
loss of generality that 
\[
\hat{g}_{k}^{j}(z_{k})\rightarrow \zeta ^{(j)}\in \{-(n-1),\ldots ,0\} 
\]
with $\zeta ^{(j+1)}=\zeta ^{(j)}+1$ whenever $\zeta ^{(j)}\neq 0$. It
follows that some $\zeta ^{(j)}=0$, whence $\hat{g}_{k}^{j}(z_{k})=O(\gamma
_{k}^{-1})$ as $\hat{g}_{k}^{j+1}(z_{k})$ is bounded away from $1$.
Calculating the derivative 
\[
\hat{g}_{k}^{\prime }(z)=\frac{1}{\alpha _{k}}-\frac{\alpha _{k}}{(\alpha
_{k}\gamma _{k}-1)z^{2}} 
\]
we see that $\hat{g}_{k}^{\prime }(\hat{g}_{k}^{j}(z_{k}))\rightarrow \infty 
$ when $\zeta ^{(j)}=0$, while $\hat{g}_{k}^{\prime }(\hat{g}%
_{k}^{j}(z_{k}))\rightarrow 1$ otherwise, and we conclude that $\rho
_{k}\rightarrow \infty .\;\;\Box $

\bigskip Assume that $\alpha _{k}=\omega [1+\tau \sqrt{\epsilon _{k}}]+o(%
\sqrt{\epsilon _{k}})$, hence $\beta _{k}=\bar{\omega}[1-\tau \sqrt{\epsilon
_{k}}]+o(\sqrt{\epsilon _{k}})$, where $\omega \neq 1$ is a primitive $q$-th
root of unity and $\tau \in {\Bbb C}$. For $z\in {\Bbb C}^{*}$ and $0\leq
j<q $, it follows from (\ref{foverz}) that 
\[
f_{k}^{q}\left( \bar{\omega}^{j}[1+z\sqrt{\epsilon _{k}}]\right) =\bar{\omega%
}^{j}\left[ 1+\left( z+q\tau +\frac{1}{z+j\tau }\right) \sqrt{\epsilon _{k}}%
\right] )+o(\sqrt{\epsilon _{k}}) 
\]
whence 
\begin{equation}
\psi _{k,(j)}^{-1}\circ f_{k}^{q}\circ \psi _{k,(j)}(z)\rightarrow z+q\tau +%
\frac{1}{z+j\tau }=G_{q\tau }(z+j\tau )-j\tau  \label{fkq}
\end{equation}
locally uniformly on ${\Bbb C}^{*}$, where 
\[
\psi _{k,(j)}(z)=\bar{\omega}^{j}(1+z\sqrt{\epsilon _{k}}). 
\]
Similarly, if $\alpha _{k}\rightarrow \omega $ but $\frac{\alpha _{k}-\omega 
}{\sqrt{\epsilon _{k}}}\rightarrow \infty $ then 
\[
\frac{f_{k}^{q}\left( \bar{\omega}^{j}(1+z\sqrt{\epsilon _{k}})\right) }{%
\alpha _{k}^{q}}=\left\{ 
\begin{array}{ll}
1+\left( z+\frac{1}{z}\right) \sqrt{\epsilon _{k}}+o(\sqrt{\epsilon _{k}}) & 
j=0 \\ 
\bar{\omega}^{j}(1+z\sqrt{\epsilon _{k}})+o(\sqrt{\epsilon _{k}}) & j\neq 0
\end{array}
\right. 
\]
and thus 
\[
\psi _{k,(j)}^{-1}\circ f_{k}^{q}\circ \psi _{k,(j)}(z)\rightarrow \infty 
\]
locally uniformly on ${\Bbb C}^{*}$. Applying (\ref{magnify}) to the case $%
j=0$, we deduce:

\begin{proposition}
\label{kalimit} Let $F_{k}=F_{\gamma _{k},\delta _{k}}$ where $\gamma
_{k}\rightarrow \infty $. Assume that $\alpha _{k}\rightarrow \omega $,
hence $\beta _{k}\rightarrow \bar{\omega}$, where $\omega \neq 1$ is a
primitive $q$-th root of unity, and assume further that $\frac{\alpha
_{k}-\omega }{\sqrt{\epsilon _{k}}}\rightarrow \omega \tau $, hence $\frac{%
\beta _{k}-\bar{\omega}}{\sqrt{\epsilon _{k}}}\rightarrow -\bar{\omega}\tau $%
, for some $\tau \in \widehat{{\Bbb C}}$. Then 
\[
F_{k}^{q}\rightarrow G_{q\tau } 
\]
locally uniformly on ${\Bbb C}^{*}.$
\end{proposition}

\medskip Recalling Lemmas \ref{calz} and \ref{bddeig}, we observe:

\begin{proposition}
\label{limitcases} Let $F_{k}=F_{\gamma _{k},\delta _{k}}$ where $\gamma
_{k}\rightarrow \infty $. Assume that $\alpha _{k}\rightarrow \omega $,
hence $\beta _{k}\rightarrow \bar{\omega}$, where $\omega \neq 1$ is a
primitive $q$-th root of unity. Assume further that $F_{k}^{q}\rightarrow
G_{T}$ for some $T\in {\Bbb C}$, and let $\langle z_{k}\rangle \rightarrow 
{\cal Z}$ be cycles of period $n>1$ and eigenvalues $\rho _{k}\rightarrow
\rho _{\infty }\in {\Bbb C}$. Then $G_{T}({\cal Z})\subseteq {\cal Z}$.
Furthermore:

\begin{itemize}
\item  If ${\cal Z}=\{\infty \}$ then $T=0$.

\item  If $0\in {\cal Z}$ then $G_{T}^{m}(\chi )=0$, whence $%
G_{T}^{m+1}(\chi )=\infty =G_{T}^{m+2}(\chi )$, for some $\chi \in \{+1,-1\}$
and $1\leq m<\frac{n}{q}$.

\item  Otherwise, ${\cal Z}=\langle \zeta \rangle \cup \{\infty \}$ where $%
\langle \zeta \rangle \subset {\Bbb C}$ is a cycle of period $m=\frac{n}{q}$
and eigenvalue $\rho _{\infty }$, or possibly a parabolic cycle of lower
period if $\rho _{\infty }=1$.
\end{itemize}
\end{proposition}

\smallskip Conversely, given an $m$-cycle $\langle \zeta \rangle $ of $G_{T}$
there exist $mq$-cycles of $F_{k}$ converging to $\langle \zeta \rangle \cup
\{\infty \}$. In particular, for $T\neq 0$ there is a unique finite fixed
point $\zeta =-\frac{1}{T}$ with eigenvalue $1-T^{2}$, hence $\langle
z_{k}\rangle \rightarrow \{\zeta ,\infty \}$ for some $q$-cycles $\langle
z_{k}\rangle $. As ${\rm mult}_{G_{T}}(\zeta )=1$, it follows from Lemma \ref
{calz} that $\langle \hat{z}_{k}\rangle \rightarrow \{0,\infty \}$ for every
convergent sequence of $q$-cycles $\langle \hat{z}_{k}\rangle \neq \langle
z_{k}\rangle $. In view of Lemma \ref{bddeig}, the eigenvalues of $\langle 
\hat{z}_{k}\rangle $ tend to $\infty $, as do those of all $\ell $-cycles
where $\ell \not{\in}\{1,q\}$ divides $q$, and thus 
\begin{equation}
\frac{1}{1-\alpha _{k}^{q}}+\frac{1}{1-\beta _{k}^{q}}\rightarrow 1-\frac{1}{%
T^{2}}  \label{abqindex}
\end{equation}
by (\ref{holindex}). On the other hand, for $T=0$ there is only the fixed
point at $\infty $, so every convergent sequence of $q$-cycles of $F_{k}$
tends to $\{\infty \}$ or $\{0,\infty \}$. The validity of (\ref{abqindex})
in this case is a particular consequence of the following:

\begin{proposition}
\label{qcycles} Let $F_{k}=F_{\gamma _{k},\delta _{k}}$ where $\gamma
_{k}\rightarrow \infty $. Assume that $\alpha _{k}\rightarrow \omega $,
hence $\beta _{k}\rightarrow \bar{\omega}$, where $\omega \neq 1$ is a
primitive $q$-th root of unity, and let $\langle z_{k}\rangle $ be cycles of
period $n>1$. If $\langle z_{k}\rangle \rightarrow \{\infty \}$ then $n=q$,
and every convergent sequence of $q$-cycles $\langle \hat{z}_{k}\rangle \neq
\langle z_{k}\rangle $ tends to $\{0,\infty \}$.
\end{proposition}

\noindent {\bf Proof: } In view of Lemma \ref{pqbddeig}, we may assume
without loss of generality that $n=mq$ for some positive integer $m$. By (%
\ref{magnify}), it is enough to show that for $r$ and $k$ sufficiently large
at most one $mq$-cycle of $f_{k}$ lies completely inside 
\[
V_{k}^{r}=\widehat{{\Bbb C}}-\bigcup_{j=1}^{q-1}\overline{D}_{k,(j)}^{r} 
\]
where $D_{k,(j)}^{r}=\{z\in {\Bbb C}:\;|z-\bar{\omega}^{j}|<r\sqrt{|\epsilon
_{k}|}\}.$ It follows from (\ref{foverz}) that $f_{k}^{-mq}(\infty )\cap 
\overline{V}_{k}^{r}=\{\infty \}$ for large $r$ and $k$, and thus all of the 
$2^{mq}-1$ finite poles of $f_{k}^{mq}$ lie in $%
\bigcup_{j=0}^{q-1}D_{k,(j)}^{r}$. Consequently, 
\[
\sum_{z=f_{k}^{mq}(z)\in V_{k}^{r}}{\rm mult}_{f_{k}^{mq}}(z)\;\;=\;%
\;2^{mq}+1\;-\;\sum_{z=f_{k}^{mq}(z)\in \hat{{\bf C}}-\overline{V}_{k}^{r}}%
{\rm mult}_{f_{k}^{mq}}(z) 
\]
provided that $f_{k}^{mq}$ has no fixed points on $\partial V_{k}^{r}$,
whence 
\[
\sum_{z=f_{k}^{mq}(z)\in V_{k}^{r}}{\rm mult}_{f_{k}^{mq}}(z)=2\;-\;%
\sum_{j=0}^{q-1}\frac{1}{2\pi i}\int_{\partial D_{k,(j)}^{r}}\frac{%
1-(f_{k}^{mq})^{\prime }(z)}{z-f_{k}^{mq}(z)}\,dz 
\]
by the Argument Principle.

Observe that $G_{0}(z)=z+\frac{1}{z}$ has a fixed point of multiplicity $3$
at $\infty $, and thus $G_{0}^{m}$ has $2^{mq}-2$ finite fixed points and $%
2^{mq}-1$ finite poles. It follows as above that 
\[
\frac{1}{2\pi i}\int_{|z|=r}\frac{1-(G_{0}^{m})^{\prime }(z)}{z-G_{0}^{m}(z)}%
\,dz\;=\;-1 
\]
so long as $r>\max \{|z|:\;z\in {\Bbb C}\mbox{ and }z-G_{0}^{m}(z)\in
\{0,\infty \}\,\}$. In view of (\ref{fkq}), 
\begin{eqnarray*}
\frac{1}{2\pi i}\int_{\partial D_{k,(j)}^{r}}\frac{1-(f_{k}^{mq})^{\prime
}(z)}{z-f_{k}^{mq}(z)}\,dz &=&\frac{1}{2\pi i}\int_{|z|=r}\frac{1-\left(
\psi _{k,(j)}^{-1}\circ f_{k}^{mq}\circ \psi _{k,(j)}\right) ^{\prime }(z)}{%
z-\left( \psi _{k,(j)}^{-1}\circ f_{k}^{mq}\circ \psi _{k,(j)}\right) (z)}%
\,dz \\
&=&\frac{1}{2\pi i}\int_{|z|=r}\frac{1-(G_{0}^{m})^{\prime }(z)}{%
z-G_{0}^{m}(z)}\,dz
\end{eqnarray*}
when $k$ is sufficiently large, and thus 
\[
\sum_{z=f_{k}^{mq}(z)\in V_{k}^{r}}{\rm mult}_{f_{k}^{mq}}(z)=q+2. 
\]
We deduce that $f_{k}^{mq}(z)=z\in V_{k}^{r}$ implies $f_{k}^{q}(z)=z$ for
large $r$ and $k$ depending only on $m$. If $\alpha _{k}^{q}=1$ then ${\rm %
mult}_{f_{k}^{mq}}(0)=q+1$ and ${\rm mult}_{f_{k}^{mq}}(\infty )=1$, while
if $\beta _{k}^{q}=1$ then ${\rm mult}_{f_{k}^{mq}}(0)=1$ and ${\rm mult}%
_{f_{k}^{mq}}(\infty )=q+1$; in these cases $f_{k}^{mq}$ has no other fixed
points in $V_{k}^{r}$. Otherwise, 
\[
{\rm mult}_{f_{k}^{mq}}(0)=1={\rm mult}_{f_{k}^{mq}}(\infty ) 
\]
and it follows from (\ref{foverz}) that the remaining $q$ fixed points of $%
f_{k}^{q}$ in $V_{k}^{r}$ constitute a $q$-cycle of $f_{k}$. $\Box $

\medskip

In view of Fatou's Theorem, the second assertion in Proposition \ref
{limitcases} is sharpened by:

\begin{proposition}
\label{bshrink} Let $F_{k}=F_{\gamma _{k},\delta _{k}}$ where $\gamma
_{k}\rightarrow \infty $, and let $z_{k}$ be attracting points of period $%
n>1 $ with immediate basins $B_{k}$. If $z_{k}\rightarrow 0$ then $%
B_{k}\rightarrow 0$.
\end{proposition}

\noindent {\bf Proof: } In view of Proposition {\bf \ref{no01} }we may
assume without loss of generality that $\alpha _{k}=\omega +O(\sqrt{\epsilon
_{k}})$ where $\omega \neq 1$ is a root of unity. If $k$ is large then $%
z_{k}\in {\Bbb D}$, so for $j\geq 0$ there are unique components $%
W_{k}^{j}\ni z_{k}$ of $F_{k}^{-nj}({\Bbb D}).$ We claim first that $%
W_{k}^{1}\rightarrow 0$; otherwise, as $W_{k}^{1}$ is connected there exist $%
k_{\ell }\rightarrow \infty $ and $w_{k_{\ell }}\in W_{k_{\ell }}^{1}$ with $%
\sqrt{\epsilon _{k_{\ell }}}=o(w_{k_{\ell }})$, contradicting Lemma \ref
{scaledown}. In particular, $W_{k}^{1}\subset {\Bbb D}$ and thus $%
W_{k}^{j+1}\subset W_{k}^{j}$ for $j\geq 0$ and sufficiently large $k$. We
denote $W_{k}^{\infty }$ the component of $z_{k}$ in the interior of $%
\bigcap_{j=0}^{\infty }W_{k}^{j}$ and contend that $W_{k}^{\infty }=B_{k}.$
By definition, if $\zeta \in B_{k}$ there exists open $U\ni \zeta $ such
that $F_{k}^{nj}(U)\subset {\Bbb D}$ when $j$ is large, while if $\zeta \in
\partial W_{k}^{\infty }$ there exist $\zeta _{j}\rightarrow \zeta $ with $%
F_{k}^{nj}(\zeta _{j})\in \partial {\Bbb D}$. Thus, $B_{k}\cap \partial
W_{k}^{\infty }=\emptyset $ and as $z_{k}\in B_{k}$ it follows that $%
B_{k}\subseteq W_{k}^{\infty }$; conversely, $W_{k}^{\infty }\subseteq B_{k}$
as $F_{k}^{nj}$ is bounded, hence normal, on $W_{k}^{\infty }$. $\Box $

\section{Precompactness}

\label{precom}

\medskip Recall that a rational map is {\em hyperbolic} if and only if the
orbit of every critical point tends to some attracting cycle. As discussed
in \cite{milnor2,rees2}, there are four configurations for quadratics:

\begin{description}
\item[{\bf B }]  Both critical points lie in the immediate basin of the same
attracting cycle, but in different components.

\item[{\bf C }]  Both critical points lie in the basin of the same
attracting cycle, but only one lies in the immediate basin.

\item[{\bf D }]  The critical points lie in the immediate basins of distinct
attracting cycles.

\item[{\bf E }]  Both critical points lie in the same component of the
immediate basin of an attracting fixed point.
\end{description}

There is in fact a unique hyperbolic component of type E consisting of maps
with totally disconnected Julia set. This component is unbounded; see \cite
{milnor2} for details. Our main result is that components of type D are
bounded, so long as neither attractor is a fixed point:

\begin{theorem}
\label{main} Let $g_{k}$ be quadratic rational maps, each having distinct
nonrepelling cycles of periods $n^{\pm }>1$. Then $[g_{k}]$ is bounded in $%
{\bf Rat}_{2}$.
\end{theorem}

\noindent {\bf Proof: } It follows from the proof of the Fatou-Shishikura
Inequality that we lose no generality in assuming that these cycles are
attracting \cite{shish}. Suppose to the contrary that $[g_{k}]$ is unbounded
in ${\bf Rat}_{2}$. By Lemma \ref{unbdd} we may, passing to a subsequence if
necessary, choose representatives $F_{k}=F_{\gamma _{k},\delta _{k}}\in
[g_{k}]$ where $\gamma _{k}\rightarrow \infty $. Let $\langle z_{k}^{\pm
}\rangle $ be the corresponding $n^{\pm }$-cycles of $F_{k}$. Without loss
of generality, $\langle z_{k}^{\pm }\rangle \rightarrow {\cal Z}^{\pm
}\subset \widehat{{\Bbb C}}$; it follows from Propositions \ref{no01} and 
\ref{kalimit} that we may also assume that $\alpha _{k}\rightarrow e^{2\pi
ip/q}$ and $F_{k}^{q}\rightarrow G_{T}$ for some $q\geq 2$ and $T\in {\Bbb C}
$. In view of Fatou's Theorem we may label the critical points so that $\pm
1 $ lies in the immediate basin of $\langle z_{k}^{\pm }\rangle $.

By Proposition \ref{limitcases}, if ${\cal Z}^{\pm }=\{\infty \}$ then $T=0$
and thus also ${\cal Z}^{\mp }=\{\infty \}$ by Lemma \ref{fsineq}, in
contradiction to Proposition \ref{qcycles}. Thus, ${\cal Z}^{\pm }=\langle
\zeta ^{\pm }\rangle \cup \{\infty \}$ for some nonrepelling cycle $\langle
\zeta ^{\pm }\rangle \subset {\Bbb C}$ of $G_{T}$, or else $G_{T}^{m_{\pm
}}(\pm 1)=0$ for some $m_{\pm }\geq 1$ as a consequence of Proposition \ref
{bshrink}. Applying Lemma \ref{fsineq}, we deduce that ${\cal Z}%
^{+}-\{\infty \}=\langle \zeta \rangle ={\cal Z}^{-}-\{\infty \}$ for some
cycle $\langle \zeta \rangle $. It follows from Lemma \ref{paratt} that $%
\langle \zeta \rangle $ is parabolic-attracting or parabolic-indifferent,
once again contradicting Lemma \ref{fsineq}. $\Box $

\bigskip The same considerations apply when there is one nonrepelling cycle
along with a preperiodic critical point:

\begin{theorem}
\label{slice} Let $g_{k}$ be quadratic rational maps with nonrepelling $n$%
-cycles $\langle z_{k}\rangle $ where $n>1$. Assume further that $%
g_{k}^{\ell }(\chi _{k})\in \langle \hat{z}_{k}\rangle $ for some $\ell >0$,
critical points $\chi _{k}$, and $\hat{n}$-cycles $\langle \hat{z}%
_{k}\rangle $. Then $[g_{k}]$ is bounded in ${\bf Rat}_{2}$.
\end{theorem}

The exceptional type D components are known to be unbounded; see Lemma \ref
{caryoc} below. Many, though not all type D maps arise as {\em matings} of
pairs of hyperbolic quadratic polynomials. In this construction, the
filled-in Julia sets are glued back-to-back along complex-conjugate prime
ends; see \cite{2mate,tanlei} for further details. It is tempting to
speculate that our arguments could be refined to establish precompactness
for large portions of the mating locus, but our results in this direction
are rather limited at present

\begin{figure}[tbp]
\centerline{\psfig{figure=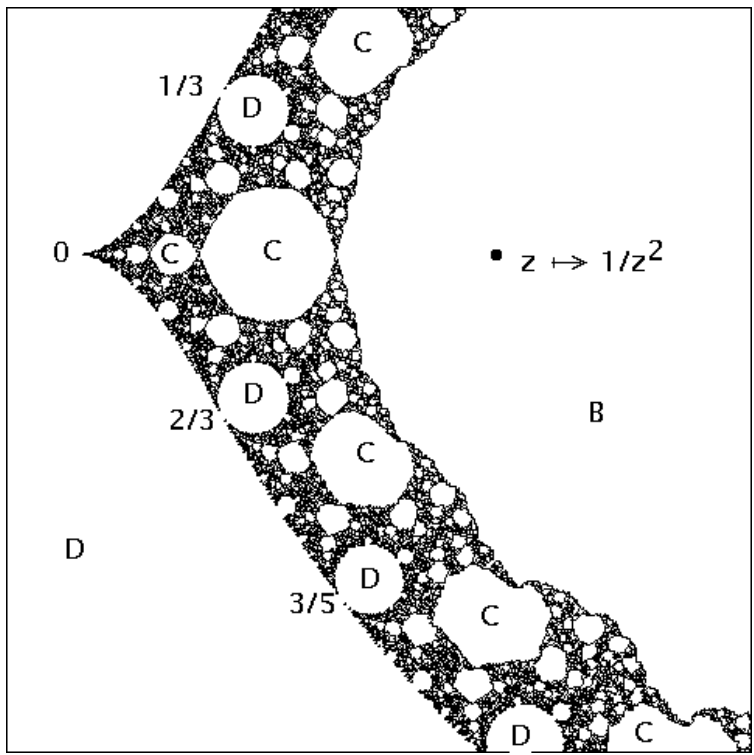,width=.6\hsize}}
\caption{Bifurcation locus in ${\rm Per}_2(0)$.}
\label{per20}
\end{figure}

Examination of Figure \ref{per20} suggests that the type C components are
all bounded. This would follow immediately from our arguments if it could be
shown in this case that $F_{k}^{-1}({\cal B}_{k})-{\cal B}_{k}\rightarrow 0$
where ${\cal B}_{k}$ is the immediate basin of the unique attracting cycle.
There are evidently many unbounded type B components. Makienko \cite
{makienko} has obtained a degree-independent sufficient condition for
unboundedness, loosely speaking the existence of a family of closed
Poincar\'{e} geodesics on the basin quotient with lifts linking to separate
the Julia set; see also \cite{pilgrim}. On the other hand, there are type B
maps which do not admit such a family: Pilgrim cites the example $g(z)=\frac{%
i\sqrt{3}}{2}(z+\frac{1}{z})$ and describes its Julia set as an {\em %
almost-Sierpinski carpet}. Such maps presumably lie in bounded hyperbolic
components.

A good deal of what is known about hyperbolic quadratic rational maps - that
Fatou components are usually Jordan domains \cite{pilgrim2}, that
polynomials can be mated if and only if they do not lie in conjugate limbs
of the Mandelbrot set \cite{tanlei}, that mating is discontinuous due to the
existence of type D hyperbolic components whose closures are not
homeomorphic to $\overline{{\Bbb D}}\times \overline{{\Bbb D}}$ \cite{2mate}%
, that moduli space is isomorphic to ${\Bbb C}^{2}$ - is valid with minor
changes for higher degree {\em bicritical maps} possessing two maximally
degenerate critical points. Much of the discussion here extends similarly,
and Milnor has recently generalized Lemma \ref{unbdd} to this larger
setting: if $[g_{k}]$ is unbounded then the eigenvalues of all but at most
two fixed points tend to infinity \cite{milcom}. However, it is not
immediately apparent how best to adapt the brute-force calculations of
Section \ref{estsec}, or better yet, how to replace them with a more
conceptual approach applicable to other degenerating families.

\section{Intersection Theory}

\label{conclusion}

The results above yield preliminary information about the intersection
theory at infinity of dynamically defined curves in moduli space. Milnor's
isomorphism $j:{\bf Rat}_{2}\rightarrow {\Bbb C}^{2}$ induces a natural
compactification $\widehat{{\bf Rat}}_{2}\cong {\Bbb P}^{2}$. Following the
discussion in \cite{milnor1} we identify the line at infinity ${\cal L}$
with the set of unordered triples $[\alpha ,\alpha ^{-1},\infty ]$ where $%
\alpha \in \widehat{{\Bbb C}}$, so that $\alpha +\alpha ^{-1}$ is the
limiting ratio of $\frac{Y}{X-2}$ in the coordinates of Section \ref
{paramspace}; see \cite{silver} for a treatment in language of geometric
invariant theory. With this convention, an unbounded sequence $[g_{k}]\in 
{\bf Rat}_{2}$ converges to the ideal point $[\alpha ,\alpha ^{-1},\infty ]$
if and only if $[\alpha _{k},\beta _{k},\gamma _{k}]\rightarrow [\alpha
,\alpha ^{-1},\infty ]$, where $\alpha _{k},\beta _{k},\gamma _{k}$ are the
eigenvalues of the fixed points of $g_{k}$. Note that the degeneration
described in Proposition \ref{kalimit} takes place in a parameter space
where 
\[
\infty _{p/q}=[e^{2\pi ip/q},e^{-2\pi ip/q},\infty ]=\infty _{(q-p)/q} 
\]
has been blown up and replaced by a 2-fold branched cover of the line ${\rm %
Per}_{1}(1).$

Recall that a {\em curve} $C$ in ${\Bbb P}^{2}$ may be defined as an
equivalence class of homogeneous polynomials in ${\Bbb C}[W,X,Y]$, where $%
H\sim \tilde{H}$ when $H=\lambda \tilde{H}$ for some $\lambda \in {\Bbb C}%
^{*}$. We write $C=\underline{H}$ and $\deg \;C=\deg \;H$; a point $P\in 
{\Bbb P}^{2}$ with homogeneous coordinates $[w:x:y]$ {\em belongs to} $C=%
\underline{H}$ if and only if $H(w,x,y)=0$, and we write $P\in C$. An {\em %
algebraic family} of degree $d$ curves parametrized by a variety $\Lambda $
is a regular map $\Lambda \rightarrow {\cal C}_{d}$, where the set ${\cal C}%
_{d}$ of all degree $d$ curves is naturally regarded as the projective space 
${\Bbb P}^{\frac{d(d+3)}{2}}$.

If $H$ has no nontrivial factors then $C=\underline{H}$ is said to be {\em %
irreducible}. An irreducible curve $\hat{C}=\underline{\hat{H}}$ with $\hat{H%
}|H$ is a {\em component} of $C$, and curves $C_{1}$ and $C_{2}$ with no
common component are said to {\em intersect properly}. Notice that $C=%
\underline{H}$ intersects ${\cal L}$ properly if and only if $\deg
\;H(1,X,Y)=\deg \;H$. Curves $C_{1},C_{2}$ which intersect properly have
finitely many points in common, and each such point can be assigned an
appropriate {\em intersection multiplicity} ${\cal I}_{C_{1},C_{2}}(P)>0$;
the {\em intersection cycle} is the formal sum 
\[
C_{1}\bullet C_{2}=\sum_{P\in C_{1}\cap C_{2}}{\cal I}_{C_{1},C_{2}}(P)\cdot
P\;. 
\]
Bezout's Theorem asserts that the total intersection multiplicity is the
product of the degrees $d_{i}=\deg \;C_{i}$, whence $C_{1}\bullet C_{2}$ may
be regarded as an element of the symmetric product 
\[
{\cal S}_{d_{1}d_{2}}\;=\;{\rm Sym}^{d_{1}d_{2}}({\Bbb P}^{2}). 
\]
Moreover, $(C_{1},C_{2})\leadsto C_{1}\bullet C_{2}$ yields a regular map 
\[
{\cal C}_{d_{1}}\times {\cal C}_{d_{2}}\,-\,{\cal E}_{d_{1},d_{2}}%
\rightarrow {\cal S}_{d_{1},d_{2}} 
\]
where ${\cal E}_{d_{1},d_{2}}$ is the set of pairs of curves with a common
component; see \cite{algcurve} for further details. The {\em intersection
cycle at infinity} is 
\[
C_{1}\bullet _{\infty }C_{2}\;=\;\sum_{P\in C_{1}\cap {\cal L}\cap C_{2}}%
{\cal I}_{C_{1},C_{2}}(P)\cdot P\;. 
\]

Consider the function $n\leadsto d(n)$ defined inductively by the relation 
\[
\sum_{m|n}d(m)=2^{n-1}; 
\]
equivalently, $d(n)$ is the number of period $n$ hyperbolic components of
the Mandelbrot set $M$. Milnor has shown the following \cite{milnor1}:

\begin{lemma}
\label{perdeg} For each $n\geq 1$ there is a algebraic family of curves 
\[
{\Bbb C}\ni \rho \leadsto {\rm Per}_{n}(\rho )\in {\cal C}_{d(n)} 
\]
uniquely determined by the condition that $[g]\in {\rm Per}_{n}(\rho )$ for $%
\rho \neq 1$ if and only if $g$ has an $n$-cycle with eigenvalue $\rho $.
The curves ${\rm Per}_{n}(1)$ are reducible for $n>1$, indeed 
\[
{\rm Per}_{n}(1)={\rm Per}_{n}^{\#}(1)\,\cup \,\bigcup_{1<q|n,\,(p,q)=1}{\rm %
Per}_{\frac{n}{q}}(e^{2\pi ip/q}) 
\]
where the generic $[g]\in {\rm Per}_{n}^{\#}(1)$ has an $n$-cycle of
eigenvalue 1.
\end{lemma}

\medskip Here are the defining polynomials for $n=1,2,3$: 
\[
\begin{array}{ll}
{\rm Per}_{1}(\rho ):\;\; & \rho ^{3}W-\rho ^{2}X+\rho Y-X+2W \\ 
{\rm Per}_{2}(\rho ):\;\; & \rho W-2X-Y \\ 
{\rm Per}_{3}(\rho ):\;\; & \rho ^{2}W^{3}-\rho
[WX(2X+Y)+3W^{2}X+2W^{3}]+(X+Y)^{2}(2X+Y) \\ 
& \mbox{}-WX(X+2Y)+12W^{2}X+28W^{3}
\end{array}
\]
Notice that 
\[
{\rm Per}_{1}(\rho )\bullet {\cal L}=[\rho ,\rho ^{-1},\infty ]. 
\]
For $n>1$, it follows from Proposition \ref{no01} that ${\rm Per}_{n}(\rho
)\bullet {\cal L}$ consists of points of the form $\infty _{p/q}$ where $%
1\leq p<q\leq n$; it is easily verified that 
\begin{equation}
{\rm Per}_{2}(\rho )\bullet {\cal L}=\infty _{1/2}  \label{per2}
\end{equation}
and 
\begin{equation}
{\rm Per}_{3}(\rho )\bullet {\cal L}=\infty _{1/2}+2\cdot \infty _{1/3}
\label{per3}
\end{equation}
for every $\rho \in {\Bbb C}$. The degeneration described in Proposition \ref
{kalimit} takes place in a parameter space where $\infty _{p/q}$ has been
blown up and replaced by a 2-fold branched cover of the line ${\rm Per}%
_{1}(1)$.

Recall that the $p/q-${\em limb} of $M$ is the set $L_{p/q}$ of all
parameter values for which $P_{c}(z)=z^{2}+c$ has a fixed point of
combinatorial rotation number $p/q$; see for example \cite{hubbard}. Given $%
\alpha \in {\Bbb D}$, it follows from standard deformation considerations 
\cite{mss} that there is a unique class $[P_{c,(\alpha )}]\in {\rm Per}%
_{1}(\alpha )$ consisting of maps which are quasiconformally conjugate to $%
P_{c}$ on a neighborhood of the filled-in Julia set $K(P_{c})$ through
conjugacies with vanishing dilatation on $K(P_{c})$. Petersen \cite
{milnor1,carsten} showed the following by a modulus estimate similar to that
in the proof of Yoccoz Inequality:

\begin{lemma}
\label{caryoc} Let $P_{c}(z)=z^{2}+c$ where $c\in M$. If $\alpha _{k}\in 
{\Bbb D}$ converges nontangentially to $e^{-2\pi ip/q}\neq 1$ then $%
[P_{c_{k},(\alpha _{k})}]\rightarrow \infty _{p/q}\in {\cal L}$ uniformly
for $c_{k}\in L_{p/q}.$
\end{lemma}

\medskip Each $c\in L_{p/q}$ determines an arc $\left\{ [P_{c,(\alpha
)}]:\;e^{2\pi ip/q}\alpha \in [0,1)\right\} $ with an endpoint at $[P_{c}]$.
In view of Lemma \ref{caryoc} the other endpoint is $\infty _{p/q}$; in
particular, each of the period $n$ components in $L_{p/q}$ yields a branch
of ${\rm Per}_{n}(\rho )$ at $\infty _{p/q}$, at least for $\rho \in 
\overline{{\Bbb D}}{\bf .}$ As distinct $P_{c}$ lie on disjoint arcs, it
follows from Lemma \nolinebreak \ref{perdeg} and Bezout's Theorem that there
are $d(n)$ such branches in total. Thus, 
\begin{equation}
{\rm Per}_{n}(\rho )\bullet {\cal L}\;=\;\sum_{1<q|n,\,(p,q)=1}d_{p/q}(n)%
\cdot \infty _{p/q}  \label{pern0l}
\end{equation}
where $d_{p/q}(n)$ is the cardinality of 
\[
L_{p/q}(n)=\{W:\;W\subset L_{p/q}\mbox{ is a period $n$ component 
of }M\}; 
\]
as $\rho \leadsto {\rm Per}_{n}(\rho )$ is continuous, (\ref{pern0l}) holds
for every $\rho \in {\Bbb C}$, in accordance with Proposition \ref{no01}.
Moreover, it follows conversely that $[P_{c_{k},(\alpha _{k})]}$ is bounded
in ${\bf Rat}_{2}$ if $c_{k}\in W\subset L_{p/q}{\bf \ }$and $\alpha _{k}\in 
{\Bbb D}$ is bounded away from $e^{-2\pi ip/q}$. Indeed,\ if $%
[P_{c_{k},(\alpha _{k})]}\in {\rm Per}_{n}(\rho _{k})$ is unbounded then
necessarily $[P_{c_{k},(\alpha _{k})]}\rightarrow \infty _{p/q}$, hence $%
\alpha _{k}\rightarrow e^{\pm 2\pi ip/q}$ after passage to a subsequence;
the sign is determined by the fact the $d_{p/q}(n)$ nearby points of ${\rm %
Per}_{1}(\alpha _{k})\bullet $ ${\rm Per}_{n}(\rho _{k})$ are all
deformations of polynomials $P_{c}$ with $c\in L_{p/q}$ rather than $%
L_{-p/q} $ .

\begin{figure}[tbp]
\centerline{\psfig{figure=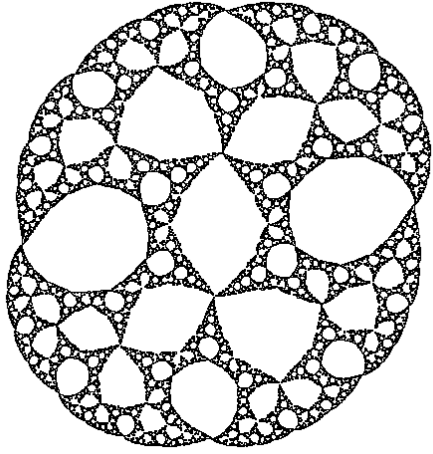,width=.4\hsize}}
\caption{$J(f)$ for $f\in{\cal RAT}_2$ with critical points of periods 2 and
3.}
\label{matefig}
\end{figure}

As shown in \cite{2mate}, the intersection of ${\rm Per}_{n^{+}}(\rho ^{+})$
and ${\rm Per}_{n^{-}}(\rho ^{-})$ is generically proper. It is somewhat
surprising that there are nontrivial exceptions: as observed in \cite
{milnor1}, it follows from (\ref{holindex}) that 
\begin{equation}
{\rm Per}_{2}(-3)={\rm Per}_{3}^{\#}(1).  \label{per2per3}
\end{equation}
In conjunction with the explicit expressions following Lemma \ref{perdeg},
this coincidence yields a short independent proof of Theorem \ref{main} in
the special case $(n^{+},n^{-})=(2,3)$. There are exactly two such
hyperbolic components, a complex-conjugate pair obtained by mating the
unique period 2 component in $M$ with the period 3 components in $L_{1/3}$
and $L_{2/3}$; see Figure~\ref{matefig}.

Recall that a quadratic rational map has precisely two $3$-cycles counting
multiplicity, whence 
\[
{\rm Per}_{2}(-3)\bullet {\rm Per}_{3}(\rho ^{-})={\rm Per}%
_{3}^{\#}(1)\bullet {\rm Per}_{3}(\rho ^{-})=3\cdot \infty _{1/2} 
\]
for $\rho ^{-}\neq 1$ by (\ref{per2}), (\ref{per3}), (\ref{per2per3}) and
Bezout's Theorem; thus, ${\rm Per}_{2}(-3)$ and ${\rm Per}_{3}(\rho ^{-})$
are tangent at $\infty _{1/2}$. In view of the transversality of distinct
lines ${\rm Per}_{2}(\rho ^{+})$, the curves ${\rm Per}_{2}(\rho ^{+})$ and $%
{\rm Per}_{3}(\rho ^{-})$ are transverse at $\infty _{1/2}$ provided that
they intersect properly. Consequently, 
\[
{\rm Per}_{2}(\rho ^{+})\bullet _{\infty }{\rm Per}_{3}(\rho ^{-})=\infty
_{1/2} 
\]
for $(\rho ^{+},\rho ^{-})\neq (-3,1)$; in particular, for $(\rho ^{+},\rho
^{-})\in \overline{{\Bbb D}}\times \overline{{\Bbb D}}$ the points in 
\[
{\rm Per}_{2}(\rho ^{+})\bullet {\rm Per}_{3}(\rho ^{-})\;-\;{\rm Per}%
_{2}(\rho ^{+})\bullet _{\infty }{\rm Per}_{3}(\rho ^{-}) 
\]
are uniformly bounded away from ${\cal L}$.

Conversely, it follows from Theorem \ref{main} and the remarks after (\ref
{pern0l}) that 
\[
{\rm Per}_{n^{+}}(\rho ^{+})\bullet _{\infty }{\rm Per}_{n^{-}}(\rho
^{-})\;=\;\sum_{\stackrel{1\leq p<q\leq \min (n^{+},n^{-})}{(p,q)=1}}{\cal I}%
_{p/q}\cdot \infty _{p/q} 
\]
independent of $\rho ^{\pm }$, for $n^{\pm }>1$ and generic $(\rho ^{+},\rho
^{-})\in {\Bbb C}^{2}$. Heuristic considerations supported by calculations
in \cite{stim} suggest that 
\[
{\cal I}_{p/q}\;=\;\sum_{(W^{+},W^{-})\in L_{p/q}(n^{+})\times
L_{p/q}(n^{-})}\iota (W^{+},W^{-}) 
\]
where $\iota (W^{+},W^{-})$ measures the mutual combinatorial depth of $%
W^{\pm }$ in $M$. The language of {\em internal addresses} \cite{dierk} may
be useful in the formulation and proof of this assertion.

\bigskip\noindent{\bf Address:}

\smallskip\noindent Institute for Mathematical Sciences

\noindent Department of Mathematics

\noindent State University of New York at Stony Brook

\noindent Stony Brook, NY 11794-3651

\medskip \noindent{\bf E-mail: } adame@math.sunysb.edu

\end{document}

%% file: imsmark.tex
\def\SBIMSMark#1#2#3{
 \font\SBF=cmss10 at 10 true pt
 \font\SBI=cmssi10 at 10 true pt
 \setbox0=\hbox{\SBF Stony Brook IMS Preprint \##1}
 \setbox2=\hbox to \wd0{\hfil \SBI #2}
 \setbox4=\hbox to \wd0{\hfil \SBI #3}
 \setbox6=\hbox to \wd0{\hss
             \vbox{\hsize=\wd0 \parskip=0pt \baselineskip=10 true pt
                   \copy0 \break%
                   \copy2 \break%
                   \copy4 \break}}
 \dimen0=\ht6   \advance\dimen0 by \vsize \advance\dimen0 by 8 true pt
                \advance\dimen0 by -\pagetotal
 \dimen2=\hsize \advance\dimen2 by .25 true in
  \openin2=publishd.tex
  \ifeof2\setbox0=\hbox to 0pt{}
  \else 
     \setbox0=\hbox to 3.1 true in{
                \vbox to \ht6{\hsize=3 true in \parskip=0pt  \noindent  
                \input publishd.tex 
                \vfill}}
  \fi
  \closein2
  \ht0=0pt \dp0=0pt
 \ht6=0pt \dp6=0pt
 \setbox8=\vbox to \dimen0{\vfill \hbox to \dimen2{\copy0 \hss \copy6}}
 \ht8=0pt \dp8=0pt \wd8=0pt
 \copy8
 \message{*** Stony Brook IMS Preprint #1, #2 ***}
}